\documentclass[11pt]{amsart}
\newtheorem{theorem}{Theorem}[section]

\newtheorem{lemma}[theorem]{Lemma}
\newtheorem{proposition}[theorem]{Proposition}

\theoremstyle{definition}

\theoremstyle{remark}
\newtheorem{remark}{Remark}
\numberwithin{equation}{section}
\begin{document}

\title[Sharp Liouville theorems]
{Sharp  Liouville theorems}
\author{Salvador Villegas}
\thanks{The author has been supported by the Ministerio de Ciencia, Innovaci\'on y Universidades of Spain PGC2018-096422-B-I00 and by the Junta de Andaluc\'{\i}a  FQM-116.}
\address{Departamento de An\'{a}lisis
Matem\'{a}tico, Universidad de Granada, 18071 Granada, Spain.}
\email{svillega@ugr.es}

\begin{abstract}
Consider the equation div$(\varphi^2 \nabla \sigma)=0$ in $\mathbb{R}^N,$ where $\varphi>0$. Berestycki, Caffarelli and Nirenberg \cite{BCN} proved that if there  exists $C>0$ such that  $\int_{B_R}(\varphi \sigma)^2 \leq CR^2$ for every $R\geq 1$ then $\sigma$  is necessarily constant. In this paper we provide necessary and sufficient conditions on $0<\Psi\in C([1,\infty))$ for which this result remains true if we replace $R^2$ with $\Psi(R)$  in any dimension $N$. In the case of the convexity of $\Psi$ for large $R>1$ and $\Psi'>0$, this condition is equivalent to $\displaystyle{\int_1^\infty\frac{1}{\Psi'}=\infty}$.
\end{abstract}
\maketitle
\section{Introduction and main results}
In 1978 E. De Giorgi \cite{DG} made the following conjecture:

{\bf Conjecture.} Let $u\in C^2( \mathbb{R}^N)$ be a bounded solution of the Allen-Cahn equation $-\Delta u=u-u^3$ which is monotone in one direction (for instance $\partial u/\partial_{x_N}>0$ in $\mathbb{R}^N$). Then $u$ is a 1-dimensional function (or equivalently, all level sets $\{u = s\}$ of u are hyperplanes), at least if $N\leq 8$.

\

This conjecture was proved in 1997 for $N=2$ by Ghoussoub and Gui \cite{GG}, and in 2000 for $N=3$ by Ambrosio and Cabr\'e \cite{AC} . In dimensions $N\geq 9$, del Pino, Kowalczyk, and Wei \cite{DKW} established that the conjecture does not hold, as suggested in De Giorgi's original statement. The conjecture remains still open for dimensions  $4\leq N\leq 8$.

In the proof of  the conjecture for $N\leq 3$, it is used the following Liouville-type theorem due to H. Berestycki, L. Caffarelli and L. Nirenberg \cite{BCN}:

\begin{theorem}\label{Lio}

Let $\varphi\in L_{loc}^\infty (\mathbb{R}^N)$ be a positive function. Assume that $\sigma\in H_{loc}^1 (\mathbb{R}^N)$ satisfies $\sigma \, \mbox{div}\, (\varphi^2\nabla\sigma )\geq 0$ in $\mathbb{R}^N$ in the distributional sense. For every $R>0$, let $B_R=\{ \vert x\vert<R\}$ and assume that there exists a constant independent of $R$ such that
$$\int_{B_R}(\varphi \sigma)^2 dx\leq CR^2  \ \ \ \mbox{for every }\ R\geq 1.$$

Then $\sigma$ is constant.

\end{theorem}

To deduce the conjecture for $N\leq 3$ from this theorem, the authors made the following reasoning: if $u$ is a solution in De Giorgi's conjecture, consider the functions $\varphi:=\partial u/\partial_{x_N}>0$ and $\sigma_i:=\partial_{x_i}u/ \partial_{x_N}u$, for $i=1,...,N-1$.  Since both $\partial_{x_i}u$ and $\varphi$ solves the same linear equation $-\Delta v=(1-3u^2)v$, an easy computation shows that div$(\varphi^2\nabla\sigma_i )=0$. In dimensions $N\leq 3$ it is proved that there exists $C>0$ such that $\int_{B_R}\vert\nabla u\vert^2 dx\leq CR^2$, for every $R\geq 1$. Applying Theorem \ref{Lio} gives $\sigma_i$ is constant for every $i=1,...,N-1.$ It follows easily that $u$ is a one-dimensional function. Observe that in the previous reasoning it is only used div$(\varphi^2\nabla\sigma_i )=0$, which is a an stronger condition than $\sigma_i \, \mbox{div}(\varphi^2\nabla\sigma_i )\geq 0$.

Motivated  by the useful application of  Liouville-type theorems to these kind of problems, a natural question is to find functions $0<\Psi\in C([1,\infty ))$, for which Theorem \ref{Lio} remains true if we replace $CR^2$ with $\Psi(R)$. In this way, we make the following definitions:

\

{\bf Property (P)}. We say that a function  $0<\Psi\in C([1,\infty ))$ satisfies (P) if it has the following property: if $\varphi\in L_{loc}^\infty (\mathbb{R}^N)$ is a positive function, $\sigma\in H_{loc}^1 (\mathbb{R}^N)$ satisfies

\begin{equation}\label{ine}
\sigma \, \mbox{div}\, (\varphi^2\nabla\sigma )\geq 0 \ \ \ \mbox{in }\mathbb{R}^N
\end{equation}

\noindent in the distributional sense and 
$$\int_{B_R}(\varphi \sigma)^2 dx\leq \Psi(R)  \ \ \ \mbox{for every }\ R\geq 1,$$

\noindent then $\sigma$ is necessarily constant.

\

{\bf Property (P')}. We say that a function  $0<\Psi\in C([1,\infty ))$ satisfies (P') if it has the following property: if $\varphi\in L_{loc}^\infty (\mathbb{R}^N)$ is a positive function, $\sigma\in H_{loc}^1 (\mathbb{R}^N)$ satisfies 

\begin{equation}\label{equ}
\mbox{div}\, (\varphi^2\nabla\sigma )= 0 \ \ \ \mbox{in }\mathbb{R}^N
\end{equation}

\noindent in the distributional sense and 
$$\int_{B_R}(\varphi \sigma)^2 dx\leq \Psi(R)  \ \ \ \mbox{for every }\ R\geq 1,$$

\noindent then $\sigma$ is necessarily constant.

\

Note that, a priori, the definitions of properties (P) and (P') depend on the dimension $N$. We will show that, in fact, this is not so: if a function $0<\Psi\in C([1,\infty ))$ satisfies (P) (resp. (P')) in some dimension $N_0$, then it satisfies (P) (resp. (P')) in any dimension $N$.

It is obvious that property (P) is stronger than property (P'). In fact, in this paper we will prove that they are equivalent. 

With this notation, Theorem \ref{Lio} says that the function $CR^2$ satisfies (P) for every $C>0$. In \cite{AAC} the authors formulated the following problem: What is the optimal (maximal) exponet $\gamma_N$ such that  $CR^{\gamma_N}$ $(C>0)$ satisfies (P')?

In \cite{B} it is proved that $\gamma_N<N$ when $N\geq 3$. Also, a sharp choise in the counterexamples of \cite{GG} shows that $\gamma_N<2+2\sqrt{N-1}$ for $N\geq 7$. Recently, Moradifam \cite{M} proved that  $\gamma_N<3$ when $N\geq 4$. Finally, in a recent work \cite{YO} the author has proven that $\gamma_N=2$ for every $N\geq1$. In other words, the functions $CR^k$ do not satisfy (P') for every $k>2$ and $C>0$. On the other hand, the sharpness of the exponent $2$ for condition (P) was proved by Gazzola \cite{GA}.

Moschini \cite{MO} proved that $CR^2(1+\log R)$ satisfies (P) for every $C>0$. By a classical example \cite{PW} it is obtained that $R^2(1+\log R)^2$ does not satisfy (P) in dimension $N=2$.

All the results previously exposed are covered by the following theorems:

\begin{theorem}\label{general}
Suppose $0<\Psi\in C([1,\infty ))$. The following conditions are equivalent:
\begin{enumerate}
\item[i)] $\Psi$ satisfies (P).
\item[ii)] $\Psi$ satisfies (P').
\item[iii)]  $\displaystyle{\int_1^\infty\frac{1}{h'}=\infty}$, for every nondecreasing function $0\leq h\in C([1,\infty ))$ satisfying $h\leq\Psi$ in $[1,\infty )$.
\end{enumerate}
\end{theorem}

Note that if $0<\Psi\in C^1([1,\infty ))$ satisfy $\Psi'>0$ in $[1,+\infty )$, the we can take $h=\Psi$ in iii), obtaining that $\int_1^\infty 1/\Psi'=\infty$ is a necessary condition to have i) and ii), but not sufficient (see Remark \ref{notsuf} below). The next theorem shows that, under convexity conditions on $\Psi$, this is also a sufficient condition to obtain i) and ii).

\begin{theorem}\label{teoremon}
Suppose $0<\Psi\in C^1([1,\infty ))$ satisfy $\Psi'>0$ in $[1,+\infty )$ and $\Psi$ is convex in $[R_0,+\infty )$ for some $R_0>1$. The following conditions are equivalent:
\begin{enumerate}
\item[i)] $\Psi$ satisfies (P).
\item[ii)] $\Psi$ satisfies (P').
\item[iii')]  $\displaystyle{\int_1^\infty\frac{1}{\Psi'}=\infty} .$
\end{enumerate}
\end{theorem}

\begin{remark} For general  $0<\Psi\in C([1,\infty ))$, it is possible to prove that if $\liminf_{x\to \infty}\Psi(x)/x^2<+\infty$ then $\Psi$ satisfy (P) and (P'). Therefore, we can restrict our attention to the case $\lim_{x\to \infty}\Psi(x)/x^2=+\infty$. Thus, the condition of convexity of $\Psi$ in Theorem \ref{teoremon} seems natural and not too restrictive. 

To see that  $\liminf_{x\to \infty}\Psi(x)/x^2<+\infty$ implies $\Psi$  satisfy (P) and (P') we will apply Theorem \ref{general}. Suppose that there exists a divergent sequence $\{ R_n\}$  and a real number $C>0$ such that $\Psi(R_n)\leq CR_n^2$, $n\geq 1$ and take a  nondecreasing function $0\leq h\in C([1,\infty ))$ satisfying $h\leq\Psi$ in $[1,\infty )$. Our purpose is to obtain $\int_1^\infty 1/h'=\infty$. To this end, take an arbitrary $R>1$ and consider $n_0\in \mathbb{N}$ such that $R_n>R$ for every $n\geq n_0$. Then

$$R_n-R=\int_R^{R_n} \sqrt{h'}\frac{1}{\sqrt{h'}}\leq \left(\int_R^{R_n} h'\right)^{1/2} \left(\int_R^{R_n} \frac{1}{h'}\right)^{1/2}$$
$$=\left(h(R_n)-h(R)\right)^{1/2} \left(\int_R^{R_n} \frac{1}{h'}\right)^{1/2}\leq (\Psi (R_n))^{1/2} \left(\int_R^{R_n} \frac{1}{h'}\right)^{1/2}$$
$$\leq C^{1/2}R_n\left(\int_R^\infty \frac{1}{h'}\right)^{1/2},$$

\noindent for every $n\geq n_0$. Hence

$$\int_R^\infty \frac{1}{h'}\geq \frac{(R_n-R)^2}{CR_n^2},$$

\noindent for $n\geq n_0$. Taking limit as $n$ tends to infinity we deduce

$$\int_R^\infty \frac{1}{h'}\geq \frac{1}{C}.$$

Since $R>1$ is arbitrary we conclude $\int_1^\infty 1/h'=\infty$, which is the desired conclusion.
\end{remark}

\begin{remark}\label{notsuf}  As said before,  if $0<\Psi\in C^1([1,\infty ))$ satisfy $\Psi'>0$ in $[1,+\infty )$, it is obvious from Theorem \ref{general} that the condition $\int_1^\infty 1/\Psi'=\infty$ is necessary to have i) and ii). To show that this is not sufficient, it suffices to construct functions $\Psi, h\in C^\infty ([1,\infty ))$ satisfying $0<h<\Psi$ and $0<\Psi', h'$ in $[1,\infty )$ such that  $\int_1^\infty 1/\Psi'=\infty$ and  $\int_1^\infty 1/h'<\infty$.

To do this, for every integer $n\geq 1$ define the function $f_n:[n,n+1/2]\rightarrow \mathbb{R}$ by

$$f_n(x):=7nx+n^3, \ n\leq x\leq n+1/2.$$ 

Clearly 

$$f_n(n+1/2)=7n(n+1/2)+n^3<f_{n+1}(n+1)=7(n+1)^2+(n+1)^3, \mbox{ for every }n\geq 1.$$

Hence, there exists $0<\Psi\in C^{\infty}  ([1,+\infty )$ satisfying $\Psi'>0$ and $\Psi(x)=f_n(x)$ for every  $n\leq x\leq n+1/2$ and $n\geq 1$.

It follows that

$$\int_1^\infty \frac{1}{\Psi'}\geq \sum_{n\geq 1}\int_n^{n+1/2}\frac{1}{\Psi'}=\sum_{n\geq 1}\frac{1}{14n}=\infty .$$

Then $\int_1^\infty 1/\Psi'=\infty$. On the other hand, take an arbitrary $x\geq 1$. Then there exists an integer $n\geq 1$ such that $n\leq x<n+1$. Thus

$$\Psi(x)\geq \Psi(n)=7n^2+n^3 \geq (n+1)^3>x^3.$$

Therefore, taking $h(x)=x^3$ we have $\int_1^\infty 1/h'<\infty$, which is our claim.
\end{remark}

\section{Proof of Theorem \ref{general}.}

\noindent\textbf{Proof of Theorem \ref{general}.}

It is evident that i)$\Rightarrow$ii). Therefore we shall have established the theorem if we prove iii)$\Rightarrow$i) and ii)$\Rightarrow$iii).

{\it \textbf{Proof of iii)$\Rightarrow$i)}}

Suppose that $0<\Psi\in C([1,\infty ))$ satisfies iii) and $\int_{B_R}(\varphi \sigma)^2 \leq \Psi (R)$ for every $R\geq 1$ where $\varphi\in L_{loc}^\infty (\mathbb{R}^N)$ is a positive function and $\sigma\in H_{loc}^1 (\mathbb{R}^N)$ satisfies (\ref{ine}) in the distributional sense. Our purpose is to obtain that $\sigma$ is constant.

 If $\inf\Psi =0$ then there exists a divergent sequence $\{ R_n\}$ such that $\Psi (R_n)$ tends to $0$ as $n$ tends to $\infty$. Thus $\int_{B_{R_n}}(\varphi \sigma)^2$ also tends to $0$, which implies $\sigma=0$.

Otherwise, let $0<m:=\inf\Psi$ and consider the function 

$$h(r):=\frac{1}{2}\int_{B_r}(\varphi \sigma )^2+\frac{1}{2}m\left( 1-e^{-r}\right) \ \ \ r\geq 1.$$

Clearly, $h\leq 1/2 \Psi +1/2 m\leq \Psi$ in $[1,\infty )$ and $h$ is a positive continuous and nondecreasing function satisfying 
$$h'(r)=\frac{1}{2}\left( \int_{\vert x\vert =r}(\varphi \sigma )^2 \right)+ \frac{1}{2} m e^{-r},$$
\noindent for almost every $r>1$. From this $\int_1^\infty1/h'=\infty$. Taking into account that $1/h'(r)\leq 2e^r/m$ for almost every $r>1$ we have $1/h'\in L_{loc}^\infty ([1,\infty))$. Thus
\begin{equation}\label{siii}
\int_R^\infty \frac{1}{h'}=\infty ,\mbox{ for every } R>1
\end{equation}

Now, for arbitrary $1<R_1<R_2$  define in the ball $B_{R_2}$ the radial function $\eta$ by

$$\eta(r):=\left\{
\begin{array}{ll}
1 & \mbox{ if } 0\leq r \leq R_1 \\ \\
\frac{\displaystyle{\int_r^{R_2}}\frac{1}{h'}}{\displaystyle{\int_{R_1}^{R_2}\frac{1}{h'}}} & \mbox{ if } R_1<r\leq R_2 \\ \\
\end{array}
\right.
$$

\noindent for every $r=\vert x\vert\leq R_2$. Multiplying (\ref{ine}) by $\eta^2$ and integrating by parts in $B_{R_2}$, we obtain

$$\int_{B_{R_2}}\eta^2 \varphi^2 \vert\nabla\sigma\vert^2\leq -2\int_{B_{R_2}}\eta \varphi^2\sigma\nabla \eta\cdot\nabla \sigma $$
$$\leq 2\left(\int_{B_{R_2}}\eta^2 \varphi^2 \vert\nabla\sigma\vert^2\right)^{1/2}\left(\int_{B_{R_2}}\varphi^2\sigma^2\vert\nabla \eta\vert^2\right)^{1/2}.$$

Therefore

$$\int_{B_{R_2}}\eta^2 \varphi^2 \vert\nabla\sigma\vert^2\leq 4\int_{B_{R_2}}\varphi^2\sigma^2\vert\nabla \eta\vert^2.$$

Thus

$$\int_{B_{R_1}}\varphi^2 \vert\nabla\sigma\vert^2\leq\int_{B_{R_2}}\eta^2 \varphi^2 \vert\nabla\sigma\vert^2\leq 4\int_{B_{R_2}}\varphi^2\sigma^2\vert\nabla \eta\vert^2$$
$$=4\int_{R_1}^{R_2}\eta'(r)^2\left( \int_{\vert x\vert =r}(\varphi \sigma )^2 \right) \, dr\leq 4\int_{R_1}^{R_2}\eta'(r)^2 2h'(r)  \, dr$$
$$=\frac{8}{\displaystyle{\left( \int_{R_1}^{R_2}\frac{1}{h'}\right)^2}}\int_{R_1}^{R_2}\frac{1}{h'(r)^2}h'(r)\, dr=\frac{8}{\displaystyle{\int_{R_1}^{R_2}\frac{1}{h'}}}.$$

Fix $R_1>1$. Applying (\ref{siii}) and taking limit in the above inequality as $R_2$ tends to $\infty$ we obtain

$$\int_{B_{R_1}}\varphi^2 \vert\nabla\sigma\vert^2=0.$$

Since $R_1>1$ is arbitrary, $\sigma$ is constant, which is the desired conclusion.

{\it \textbf{Proof of ii)$\Rightarrow$iii)}} 

Suppose that iii) does not hold. That is, there exists a nondecreasing function $0\leq h\in C([1,\infty ))$ satisfying $h\leq\Psi$ in $[1,\infty )$ and  $\displaystyle{\int_1^\infty\frac{1}{h'}<\infty}$. The proof is completed by constructing a positive function $\varphi\in L_{loc}^\infty (\mathbb{R}^N)$  and a noncontant function  $\sigma\in H_{loc}^1 (\mathbb{R}^N)$ satisfying  (\ref{equ}) in the distributional sense and $\int_{B_R}(\varphi \sigma)^2 \leq \Psi (R)$ for every $R\geq 1$.

First of all, note that $0<\lim_{r\to\infty}h(r)\leq\liminf_{r\to\infty}\Psi(r)$. Since $\Psi>0$ in $[1,\infty )$, we have that $0<m:=\inf\Psi$. Consider the odd function $\mu :\mathbb{R}\rightarrow \mathbb{R}$ such that

$$\mu(r):=\left\{
\begin{array}{ll}
\displaystyle{\frac{m}{2}(1-e^{-r})} & \mbox{ if } 0\leq r \leq 1 \\ \\
\displaystyle{\frac{1}{2}\int_1^r\min\left\{ h'(s),s^2\right\} ds+  \frac{m}{2}(1-e^{-r})}  & \mbox{ if } 1<r \\ \\
\end{array}
\right.
$$

Clearly $\mu$ is continuous and increasing in $\mathbb{R}$ and satisfies, almost everywhere, that

$$\mu'(r):=\left\{
\begin{array}{ll}
\displaystyle{\frac{m}{2}e^{-\vert r\vert }} & \mbox{ if } 0\leq \vert r\vert \leq 1 \\ \\
\displaystyle{\frac{1}{2}\min\left\{ h'(\vert r\vert),r^2\right\} +  \frac{m}{2}e^{-\vert r\vert }}  & \mbox{ if } 1<\vert r\vert \\ \\
\end{array}
\right.
$$

Therefore

$$0<\frac{1}{\mu'(r)}<\frac{2}{\min\left\{ h'(r),r^2\right\} }\leq\frac{2}{h'(r)}+\frac{2}{r^2},\, \mbox{ for every } r>1.$$

Hence $1/\mu'\in L^1(1,\infty )$ and it follows immediately $1/\mu'\in L^1(\mathbb{R} )$. For this reason, taking any $0<H\in C^\infty (\mathbb{R}^{N-1})$ satisfying $\int_{\mathbb{R}^{N-1}}H^2=1/2$, we can define the functions $\varphi, \sigma : \mathbb{R}^N \rightarrow \mathbb{R}$ by

$$\varphi (x_1,...,x_{N}):=H(x_1,...,x_{N-1})\sqrt{\mu'(x_N)}\int_{x_N}^{+\infty}\frac{dr}{\mu'(r)},$$

$$\sigma(x_1,...,x_N):=\frac{1}{\displaystyle{\int_{x_N}^{+\infty}\frac{dr}{\mu'(r)}}}.$$

(If $N=1$, then define $\varphi(x)=\sqrt{\mu'(x)}\int_{x}^{+\infty}\frac{dr}{\mu'(r)}/\sqrt{2}$ and we apply the same reasonig that in the case $N>1$).

It is easy to check that

$$0<\mu'(r)\leq \frac{1}{2}r^2 +  \frac{m}{2}e^{-\vert r\vert }, \ \ \frac{1}{\mu'(r)}\leq\frac{2}{m}e^{\vert r\vert} ,\ r\in\mathbb{R}. $$

From the above it follows that $0<\varphi\in L_{loc}^\infty (\mathbb{R}^N)$ and $\vert\nabla\sigma\vert\in L_{loc}^\infty (\mathbb{R}^N)$. Thus  $\sigma\in H_{loc}^1 (\mathbb{R}^N)$. Moreover, an easy computation shows that

$$\nabla \sigma(x_1,...,x_N)=\left( 0,...,0,\frac{1}{\mu'(x_N) \displaystyle{\left( \int_{x_N}^{+\infty}\frac{dr}{\mu'(r)}\right)^2}}\right) ,$$

$$(\varphi^2 \nabla \sigma)(x_1,...,x_N)=\left( 0,...,0,H^2(x_1,...,x_{N-1})\right),$$

\noindent which implies div$(\varphi^2 \nabla \sigma)=0$ in $\mathbb{R}^N.$

Finally taking into account that $B_R\subset \mathbb{R}^{N-1}\times (-R,R)$,  we obtain for every $R\geq 1$

$$\int_{B_R}(\varphi \sigma)^2 dx=\int_{B_R} H^2(x_1,...,x_{N-1}) \mu'(x_N)dx$$
$$\leq \int_{\mathbb{R}^{N-1}}H^2  \, d(x_1,...,x_{N-1})\int_{-R}^R \mu'(r) dr=\frac{1}{2}\left(\mu(R)-\mu(-R)\right)=\mu(R)$$
$$\leq\frac{1}{2}\int_1^R h'(s)\, ds+\frac{m}{2}(1-e^{-R})\leq\frac{h(R)}{2}+\frac{m}{2}\leq \frac{\Psi(R)}{2}+\frac{\Psi(R)}{2}=\Psi(R),$$

\noindent which completes the proof. \qed

\section{Proof of Theorem \ref{teoremon}}

\begin{proposition}\label{convex}

Let $\phi \in C^1([a,b])$ a convex function satisfying $\phi'>0$ in $[a,b]$. Then

\begin{equation}\label{main}
\int_a^b\frac{1}{\phi'}\leq \int_a^b\frac{1}{g'}
\end{equation}

\noindent for every nondecreasing function $g\in C([a,b])$ satisfying $g(a)=\phi(a)$ and $g\leq \phi$ in $[a,b]$.

Moreover, equality holds if and only if $g=\phi$.

\end{proposition}

\begin{lemma}\label{bestislinear}

Let $g\in C([a,b])$ a nondecreasing function. Let $p(x)=Ax+B$, $A>0$, $B\in \mathbb{R}$ such that $g(a)=p(a)$, $g(b)\leq p(b)$. Then 

\begin{equation}\label{aaaa}
\int_a^b\frac{1}{p'}\leq \int_a^b\frac{1}{g'}\ .
\end{equation}

Moreover, equality holds if and only if $g=p$.

\end{lemma}

\noindent\textbf{Proof.}

If $\int_a^b 1/g'=\infty$ the lemma s trivial. Otherwise, applying Cauchy-Shwartz inequality we obtain

$$b-a=\int_a^b \sqrt{g'}\frac{1}{\sqrt{g'}}\leq \left(\int_a^b g'\right)^{1/2} \left(\int_a^b \frac{1}{g'}\right)^{1/2}=\left(g(b)-g(a)\right)^{1/2} \left(\int_a^b \frac{1}{g'}\right)^{1/2}.$$

Hence

$$\int_a^b \frac{1}{g'}\geq \frac{(b-a)^2}{g(b)-g(a)}\geq \frac{(b-a)^2}{p(b)-p(a)}
=\int_a^b\frac{1}{p'}.$$

On the other hand, if equality holds then all the previous inequalities become equalities. This implies that $g(b)=p(b)$ and that $\sqrt{g'}$ is a real multiple of $1/\sqrt{g'}$. That is, $g'$ is constant and, since $g(a)=p(a)$, $g(b)=p(b)$, we obtain $g=p$. \qed

\

\begin{lemma}\label{polygonal}

Let $g\in C([a,b])$ a nondecreasing function. For $1\leq i\leq m$ consider $p_i(x)=A_i x+B_i$, $A_i>0$, $B_i\in \mathbb{R}$;  such that $p_i(a)\leq g(a)$. Define

$$\overline{g}(x):=\max\left\{ g(x),p_1(x),p_2(x),...,p_m(x)\right\} , \ a\leq x\leq b.$$

Then 

\begin{equation}\label{ahi}
\int_a^b\frac{1}{\overline{g}'}\leq \int_a^b\frac{1}{g'}\ .
\end{equation}

Moreover, if $\int_a^b 1/g'<\infty$, then equality holds if and only if $g=\overline{g}$.

\end{lemma}

\noindent\textbf{Proof.}

Note that $\overline{g}$ is a nondecreasing continuous function in $[a,b]$. Therefore, the statement of the lemma has sense. 
If $\int_a^b 1/g'=\infty$ the lemma s trivial.  Hence, we will suppose in the rest of the proof that $\int_a^b 1/g'<\infty$. The proof is by induction on $m$. 

We first prove the lemma for $m=1$. To do this, consider the open set $G=\left\{x\in (a,b): p_1(x)>g(x)\right\}$. If $G=\emptyset$, then $\overline{g}=g$ and the lemma follows. Otherwise, $G$ is the countable (possible finite) disjoint union of  open intervals. That is,  $G=\cup_{n \in X} (a_n,b_n)$, where $X\subset  \mathbb{N}$ and $p_1(a_n)=g(a_n)$, $p_1(b_n)\geq g(b_n)$ for every $n\in X$. Then

$$ \int_a^b\frac{1}{g'}-\int_a^b\frac{1}{\overline{g}'}=\int_G \left( \frac{1}{g'}-\frac{1}{p_1'}\right)=\sum_{n\in X}\int_{a_n}^{b_n} \left( \frac{1}{g'}-\frac{1}{p_1'}\right) .$$

Applying Lemma \ref{bestislinear} in each interval $(a_n,b_n)$ we conclude the lemma for the case $m=1$.

We now proceed by induction. Suppose that the lemma holds for $m-1\geq 1$ and we will prove that it holds for $m$. Define 

$$h(x):=\max\left\{ g(x),p_1(x),p_2(x),...,p_{m-1}(x)\right\} , \ a\leq x\leq b.$$

By hypothesis of induction we have 

\begin{equation}\label{abc}
\int_a^b\frac{1}{h'}\leq \int_a^b\frac{1}{g'}\ .
\end{equation}

On the other hand, note that 

$$\overline{g}(x):=\max\left\{ g(x),p_1(x),p_2(x),...,p_m(x)\right\}=\max\left\{ h(x),p_m(x)\right\} , \ a\leq x\leq b.$$

It is easily seen that $h$ is a continuous nondecreasing function satisfying $p_m(a)\leq g(a)= h(a)$. Therefore applying the case of $m=1$ (which is yet proved) to functions $h(x)$ and $p_m(x)$, we obtain

\begin{equation}\label{abcd}
\int_a^b\frac{1}{\overline{g}'}\leq \int_a^b\frac{1}{h'}\ .
\end{equation}

Combining inequalities (\ref{abc}) and (\ref{abcd}) we obtain the desired inequality (\ref{ahi}). Finally, if equality holds in (\ref{ahi}), then equalities also hold in (\ref{abc}) and (\ref{abcd}). This gives $g=h=\overline{g}$ and the proof is completed. \qed

\

\

\noindent\textbf{Proof of Proposition \ref{convex}.}

We first prove (\ref{main}) in the case $g(x)<\phi(x)$ for every $x\in (a,b)$. To do this, for every positive integer $n$,  consider a partition of the interval $(a,b]$ in $2^n$ subintervals of the same length. That is

$$(a,b]=\bigcup_{k=1}^{2^n}(x_{k-1,n},x_{k,n}];  \mbox{ where } x_{k,n}=a+k\frac{b-a}{2^n}; \ \ \ 0\leq k\leq 2^n.$$

Consider now the $2^n$ lines which are tangent to the graphic of the function $y=\phi(x)$ at $x_{k,n }$, $1\leq k\leq 2^n$. That is

$$p_{k,n}(x):=\phi'(x_{k,n})(x-x_{k,n})+\phi(x_{k,n}),\ a\leq x\leq b, \ \ 1\leq k\leq 2^n.$$

Define

$$g_n(x):=\max\left\{ g(x),p_{1,n}(x),p_{2,n}(x),...,p_{2^n,n}(x)\right\} , \ a\leq x\leq b.$$

Note that the convexity of $\phi$ gives $g_n(x)\leq \phi (x)$ for every $a\leq x\leq b$, $n\geq 1$.

We claim that $g_n\to\phi$ in $L^\infty (a,b)$ as $n \to\infty$. To do this, take an arbitrary $x\in (a,b)$. Then, for fixed $n\geq 1$, there exists $1\leq k\leq 2^n$ such that $x_{k-1,n}< x\leq x_{k,n}$. Using the convexity and monotonicty of $\phi$ we deduce

$$\phi(x)\geq g_n(x)\geq p_{k,n}(x)\geq  p_{k,n}(x_{k-1,n})=\phi'(x_{k,n})(x_{k-1,n}-x_{k,n})+\phi(x_{k,n})$$
$$\geq \phi'(b)(x_{k-1,n}-x_{k,n})+\phi(x)=-\phi'(b)\frac{b-a}{2^n}+\phi(x).$$

This gives $\Vert \phi-g_n\Vert_{L^\infty (a,b)}\leq \phi'(b)\frac{b-a}{2^n}$ and the claim is proved.

Now fix $n_0>1$ and consider $a_0=a+(b-a)/2^{n_0}$ and $b_0=b-(b-a)/2^{n_0}$. Note that $a_0=x_{2^{n-n_0},n}$ and $b_0=x_{2^n-2^{n-n_0},n}$ for every $n\geq n_0$. Since $[a_0,b_0]\subset (a,b)$ and $g<\phi$ in $(a,b)$, we deduce that there exists $\varepsilon_0>0$ (depending on $n_0$) such that $g(x)<\phi (x)-\varepsilon_0$ for every $x\in[a_0,b_0]$. Using  $g_n\to\phi$ in $L^\infty (a_0,b_0)$ we can assert that there exists $n_1\geq n_0$ (depending on $\varepsilon_0$) such that $g(x)<g_n(x)$ for every $x\in[a_0,b_0]$ and $n_1\geq n_0$. Then

$$g_n(x)=\max\left\{p_{1,n}(x),p_{2,n}(x),...,p_{2^n,n}(x)\right\} , \ a_0\leq x\leq b_0, \ n\geq n_1.$$

Consider $n\geq n_1$ and $2^{n-n_0}<k\leq 2^n-2^{n-n_0}$. Take $x\in[x_{k-1,n},x_{k,n}]$. The convexity of $\phi$ yields $g_n(x)=\max\left\{p_{k-1,n}(x),p_{k,n}(x)\right\}$ and consequently $g_n'(x)\leq \phi'(x_{k,n})$. This gives 
$$\int_{x_{k-1,n}}^{x_{k,n}} \frac{1}{g'_n}\geq \frac{x_{k,n}-x_{k-1,n}}{\phi'(x_{k,n})}.$$

Therefore, applying Lemma \ref{polygonal} in the interval $[a,b]$ it follows that

$$\int_a^b \frac{1}{g'}\geq \int_a^b \frac{1}{g'_n}\geq \int_{a_0}^{b_0} \frac{1}{g'_n}=\sum_{k=2^{n-n_0}+1}^{2^n-2^{n-n_0}}\int_{x_{k-1,n}}^{x_{k,n}} \frac{1}{g'_n}\geq\sum_{k=2^{n-n_0}+1}^{2^n-2^{n-n_0}} \frac{x_{k,n}-x_{k-1,n}}{\phi'(x_{k,n})},$$

\noindent for every $n\geq n_1$. Since $1/\phi'$ is continuous in $[a_0,b_0]$ and $x_{k,n}-x_{k-1,n}=(b-a)/2^n$ we deduce that the right term of the last inequality tends to $\int_{a_0}^{b_0}1/\phi'$ as $n$ tends to $\infty$. Thus,

$$\int_a^b \frac{1}{g'}\geq \int_{a_0}^{b_0} \frac{1}{\phi'}.$$

Finally, since $n_0>1$ is arbitrary we conclude (\ref{main}) for the case $g<\phi$ in $(a,b)$.

We now turn out to the general case  $g\leq \phi$ in $(a,b)$ and we proceed to show (\ref{main}). For this purpose, consider the open set $G=\left\{x\in (a,b): \phi(x)>g(x)\right\}$. If $G=\emptyset$, then  (\ref{main}) is trivial. Otherwise, $G$ is the countable (possible finite) disjoint union of  open intervals. That is,  $G=\cup_{n \in X} (a_n,b_n)$, where $X\subset  \mathbb{N}$, $\phi (a_n)=g(a_n)$, $\phi (b_n)\geq g(b_n)$ and $\phi>g$ in $(a_n,b_n)$ for every $n\in X$. Applying the previous case in each interval $(a_n,b_n)$ we conclude

$$ \int_a^b\frac{1}{g'}-\int_a^b\frac{1}{\phi'}=\int_G \left( \frac{1}{g'}-\frac{1}{\phi'}\right)=\sum_{n\in X}\int_{a_n}^{b_n} \left( \frac{1}{g'}-\frac{1}{\phi'}\right) \geq 0.$$

It remains to prove that equality holds in (\ref{main}) if and only if $g=\phi$. To this end suppose that we have equality in (\ref{main}) for some $g$. Take an arbitrary $x_0\in[a,b]$ and consider the function

$$g_{x_0}:=\max\left\{ g(x),\phi'(x_0)(x-x_0)+\phi(x_0)\right\} , \ a\leq x\leq b.$$

Clearly $g_{x_0}$ is nondecreasing and satisfies $g\leq g_{x_0}\leq \phi$ in $[a,b]$ and $g_{x_0}(a)=g(a)=\phi(a)$. Hence

$$\int_a^b \frac{1}{g'_{x_0}}\geq \int_a^b \frac{1}{\phi'}= \int_a^b \frac{1}{g'}$$.

Applying Lemma \ref{polygonal} yields $g=g_{x_0}$ in $[a,b]$. In particular $g(x_0)=g_{x_0}(x_0)=\max\left\{ g(x_0),\phi(x_0)\right\} =\phi(x_0)$. Since $x_0\in [a,b]$ is arbitrary we conclude that $g=\phi$ in $[a,b]$ and the proposition follows. \qed

\

\noindent\textbf{Proof of Theorem \ref{teoremon}.}

Obviously, taking $h=\Psi$ in Theorem \ref{general} it follows immediately i)$\Rightarrow$ii)$\Rightarrow$iii'). 

It remains to prove iii')$\Rightarrow$i). Suppose$\displaystyle{\int_1^\infty\frac{1}{\Psi'}=\infty}$. Using again Theorem \ref{general} what is left is to show that $\displaystyle{\int_1^\infty\frac{1}{h'}=\infty}$,  for every nondecreasing function $0\leq h\in C([1,\infty ))$ satisfying $h\leq\Psi$ in $[1,\infty )$. 

To obtain a contradiction suppose that there exists a nondecreasing function $0\leq h\in C([1,\infty ))$ satisfying $h\leq\Psi$ in $[1,\infty )$ and  $\displaystyle{\int_1^\infty\frac{1}{h'}<\infty}$. We first claim that $\displaystyle{\lim_{x\to\infty}h(x)/x=+\infty}$. Conversely, suppose that there exist $M>0$ and a divergent sequence $\{ R_n\}$ such that $h(R_n)\leq MR_n$ for every positive integer $n$. Applying Cauchy-Shwartz inequality we obtain

$$R_n-1=\int_1^{R_n} \sqrt{h'}\frac{1}{\sqrt{h'}}\leq \left(\int_1^{R_n} h'\right)^{1/2} \left(\int_1^{R_n} \frac{1}{h}\right)^{1/2}$$
$$=\left(h(R_n)-h(1)\right)^{1/2} \left(\int_1^{R_n} \frac{1}{h'}\right)^{1/2}\leq \left( M R_n\right)^{1/2}\left(\int_1^\infty \frac{1}{h'}\right)^{1/2},$$

\noindent which contradicts that $\{ R_n\}$ diverges.

Consequently there exists $R_1:=\min\left\{ R\geq R_0: h(R)=\Psi'(R_0)(R-R_0)+\Psi(R_0)\right\}$.

For every $R>R_1$ define $g_R:[R_0,R]\rightarrow \mathbb{R}$ by

$$g_R(x):=\left\{
\begin{array}{ll}
\Psi'(R_0)(x-R_0)+\Psi(R_0) & \mbox{ if } R_0\leq x \leq R_1 \\ \\
h(x) & \mbox{ if } R_1<x\leq R \\ \\
\end{array}
\right.
$$

It is easily seen that $g_R\in C([R_0,R])$ is a nondecreasing function satisfying $g_R (R_0)=\Psi (R_0)$ and $g_R\leq \Psi$ in $[R_0,R]$. Then we can apply Proposition \ref{convex} in the interval $[R_0,R]$  and obtain

$$ \int_{R_0}^R\frac{1}{\Psi'}\leq \int_{R_0}^R\frac{1}{g_R'}.$$

Hence, for arbitrary $R>R_1$, we have

$$ \int_{R_1}^R\frac{1}{h'}= \int_{R_0}^R\frac{1}{g_R'}- \int_{R_0}^{R_1}\frac{1}{g_R'}\geq \left(\int_{R_0}^R\frac{1}{\Psi'}\right) -\frac{(R_1-R_0)}{\Psi'(R_0)}.$$

Since $\displaystyle{\int_{R_0}^\infty\frac{1}{\Psi'}=\infty}$, we can take limit as $R$ tends to infty, obtaining $\displaystyle{\int_{R_1}^\infty\frac{1}{h'}\geq +\infty}$. This contradicts our assumption $\displaystyle{\int_1^\infty\frac{1}{h'}<\infty}$. \qed

\end{document}